# Double Series Involving Binomial Coefficients and the Sine Integral

John M. Campbell

maxwell8@yorku.ca

By dividing hypergeometric series representations of the inverse sine by $\sin^{-1}(x)$ and integrating, new double series representations of integers and constants arise. Binomial coefficients and the sine integral are thus combined in double series.

## 1. Double Series Transformations

Infinite series for the inverse sine (and the inverse cosine) can be transformed into double series through a process based upon Catalan's triangle.

**Theorem 1.**

$$1 = \sum_{j=0}^{\infty} \sum_{k=0}^{j} \binom{2j+1}{k} \frac{(-1)^{k+j} \mathrm{Si}(\pi(j+1-k))(2j)!\,(j-k+1)}{16^j (j!)^2 (2j+1)(2j-k+2)}$$

(1.1)

*Proof.* Consider the integral $\int \frac{x^n}{\sin^{-1} x} dx$.

$$\int \frac{x^{13}}{\sin^{-1} x} dx = \frac{429\mathrm{Si}\left(2\sin^{-1}(x)\right)}{8192} - \frac{572\mathrm{Si}\left(4\sin^{-1}(x)\right)}{8192} + \frac{429\mathrm{Si}\left(6\sin^{-1}(x)\right)}{8192} - \frac{208\mathrm{Si}\left(8\sin^{-1}(x)\right)}{8192} + \frac{65\mathrm{Si}\left(10\sin^{-1}(x)\right)}{8192} - \frac{12\mathrm{Si}\left(12\sin^{-1}(x)\right)}{8192} + \frac{\mathrm{Si}\left(14\sin^{-1}(x)\right)}{8192}$$

The sequence of numerators (1, 12, 65, 208, 429, ...) corresponds to the Catalan triangle based upon $\binom{n}{k}\frac{n-2k+1}{n-k+1}$ [4].

$$\int \frac{x^{12}}{\sin^{-1} x} dx = \frac{132\mathrm{Ci}\left(\sin^{-1}(x)\right)}{4096} - \frac{297\mathrm{Ci}\left(3\sin^{-1}(x)\right)}{4096} + \frac{275\mathrm{Ci}\left(5\sin^{-1}(x)\right)}{4096} - \frac{154\mathrm{Ci}\left(7\sin^{-1}(x)\right)}{4096} + \frac{54\mathrm{Ci}\left(9\sin^{-1}(x)\right)}{4096} - \frac{11\mathrm{Ci}\left(11\sin^{-1}(x)\right)}{4096} + \frac{\mathrm{Ci}\left(13\sin^{-1}(x)\right)}{4096}$$

Again, the sequence of numerators (1, 11, 54, 154, 275, ...) corresponds to the Catalan triangle based upon $\binom{n}{k}\frac{n-2k+1}{n-k+1}$ [4]. It is easily established that:

$$\int \frac{x^{2n}}{\sin^{-1} x} dx = \frac{1}{2^{2n}} \sum_{k=0}^{n} \binom{2n}{k} \frac{2n-2k+1}{2n-k+1} \mathrm{Ci}\left((2n+1-2k)\sin^{-1}(x)\right)(-1)^{k+n}$$

$$\int \frac{x^{2n+1}}{\sin^{-1}x}dx = \frac{1}{2^{2n+1}}\sum_{k=0}^{n}\binom{2n+1}{k}\frac{2n-2k+2}{2n-k+2}\mathrm{Si}\big((2n+2-2k)\sin^{-1}(x)\big)(-1)^{k+n}$$

In [3] it is indicated that for $x^2 \leq 1$:

$$\sin^{-1}x = \sum_{j=0}^{\infty}\frac{(2j)!}{2^{2j}(j!)^2(2j+1)}x^{2j+1}$$

(1.2)

$$1 = \sum_{j=0}^{\infty}\frac{(2j)!}{2^{2j}(j!)^2(2j+1)}\frac{x^{2j+1}}{\sin^{-1}x}$$

Integrating the above infinite sum using the following equation, it follows that **Theorem 1** holds:

$$\int_0^1 \frac{x^{2n+1}}{\sin^{-1}x}dx = \frac{1}{2^{2n+1}}\sum_{k=0}^{n}\binom{2n+1}{k}\frac{2n-2k+2}{2n-k+2}\mathrm{Si}(\pi(n+1-k))(-1)^{k+n}$$

(1.3)

**Theorem 2.**

$$\frac{\pi}{2} = \sum_{k=0}^{\infty}\sum_{j=0}^{k+1}\binom{2k+3}{j}\frac{(-1)^{j+k+1}\mathrm{Si}(\pi(k+2-j))(k!)^2(k+2-j)}{(2k+1)!\,(k+1)(2k+4-j)}$$

(2.1)

*Proof*. In [3] it is indicated that for $x^2 \leq 1$:

$$(\sin^{-1}x)^2 = \sum_{k=0}^{\infty}\frac{2^{2k}(k!)^2x^{2k+2}}{(2k+1)!\,(k+1)}$$

$$x\sin^{-1}x = \sum_{k=0}^{\infty}\frac{2^{2k}(k!)^2}{(2k+1)!\,(k+1)}\frac{x^{2k+3}}{\sin^{-1}x}$$

(2.2)

$$\int_0^1 \frac{x^{2k+3}}{\sin^{-1}x}dx = \frac{1}{2^{2k+3}}\sum_{j=0}^{k+1}\binom{2k+3}{j}\frac{2k+4-2j}{2k+4-j}\mathrm{Si}(\pi(k+2-j))(-1)^{j+k+1}$$

Thus, inserting the above variation of (1.3) into (2.2), it follows that **Theorem 2** holds.

**Theorem 3.**

$$\frac{\pi}{4}=\sum_{n=1}^{\infty}\sum_{k=0}^{n}\frac{(-1)^{k+n}\mathrm{Si}(\pi(n+1-k))\binom{2n+1}{k}(n-k+1)}{n^2\binom{2n}{n}(2n-k+2)}$$

(3.1)

*Proof*. In [2] it is indicated that $(\sin^{-1}x)^2=\frac{1}{2}\sum_{n=1}^{\infty}\frac{(2x)^{2n}}{n^2\binom{2n}{n}}$.

$$2x\sin^{-1}x=\sum_{n=1}^{\infty}\frac{2^{2n}}{n^2\binom{2n}{n}}\frac{x^{2n+1}}{\sin^{-1}x}$$

Inserting (1.3) into the above infinite sum, it follows that **Theorem 3** holds.

The $\int\frac{x^{2n}}{\sin^{-1}x}dx$ equivalent of (1.3) proves to be more challenging, as Ci(0) = -∞. Cosine integral variations of the above sums (**Theorem 1** – **Theorem 3**) are surely fecund, yielding similar results.

The following relationship is easily established:

$$\int\frac{x^n}{\ln(x+1)}dx=\sum_{k=0}^{n}\binom{n}{k}\mathrm{Ei}\big((k+1)\ln(x+1)\big)(-1)^{k+n}$$

(3.2)

Thus, dividing series such as $\frac{1}{2}(\ln(1\pm x))^2=\sum_{j=1}^{\infty}\frac{(\mp 1)^{j+1}x^{j+1}}{j+1}H_j$ [3] by $\ln(1+x)$, similar results follow. Note that the function Ei represents the exponential integral defined as follows:

$$\mathrm{Ei}(z)=\sum_{k=1}^{\infty}\frac{z^k}{kk!}+\gamma+\frac{1}{2}\left(\ln(z)-\ln\left(\frac{1}{z}\right)\right)$$

## 2. A Triangular Sequence

The generalization of the integral $\int\frac{x^{2k}}{(\sin^{-1}x)^2}dx$ involves a triangular sequence. Consider the following integrals:

$$\int\frac{x^6}{(\sin^{-1}x)^2}dx=\frac{1}{64}\left(-5\mathrm{Si}(\sin^{-1}x)+27\mathrm{Si}(3\sin^{-1}x)-25\mathrm{Si}(5\sin^{-1}x)+7\mathrm{Si}(7\sin^{-1}x)-\frac{64\sqrt{1-x^2}x^6}{\sin^{-1}x}\right)$$

$$\int\frac{x^8}{(\sin^{-1}x)^2}dx=\frac{1}{256}\left(-14\mathrm{Si}(\sin^{-1}x)+84\mathrm{Si}(3\sin^{-1}x)-100\mathrm{Si}(5\sin^{-1}x)+49\mathrm{Si}(7\sin^{-1}x)-9\mathrm{Si}(9\sin^{-1}x)-\frac{256\sqrt{1-x^2}x^8}{\sin^{-1}x}\right)$$

The sequences of integers preceding the sine integral, (5, 27, 25) and (14, 84, 100, 49), correspond to the triangular sequence which is based upon $\frac{(n+1-2m)^2}{n+1-m}\binom{n}{m}$ [1]. The following generalization is easily established:

$$\int \frac{x^{2k}}{(\sin^{-1}x)^2}dx = \frac{1}{2^{2k}}\left(-\frac{2^{2k}\sqrt{1-x^2}x^{2k}}{\sin^{-1}x} + (-1)^{k+1}(2k+1)\mathrm{Si}\big((2k+1)\sin^{-1}x\big)\right.$$
$$\left. + \sum_{n=1}^{k}(-1)^n\frac{(1-2n)^2}{k+1-n}\binom{2k}{k+n}\mathrm{Si}((2n-1)\sin^{-1}x)\right)$$

Thus, new double series involving binomial coefficients and the sine integral can be established, for example, by multiplying both sides of (1.2) by $\frac{x}{(\sin^{-1}x)^2}$ and integrating.